\documentclass[twoside,reqno]{amsart}

\usepackage{amsmath}
\usepackage{amssymb}
\usepackage{amsthm}

\newcommand{\ino}{\int_0 ^{\infty}}
\newcommand{\il}{(1+\lambda t)}
\newcommand{\iol}{\frac{1}{\lambda}}

\newtheorem{theorem}{Theorem}[section]

\theoremstyle{definition}

\numberwithin{equation}{section}

\begin{document}

\title[Degenerate Euler zeta function]{Degenerate Euler zeta function}

\author{Taekyun Kim}
\address{${}^1$ Department of Mathematics, Kwangwoon University, Seoul 01897, Republic of Korea.}
\email{tkkim@kw.ac.kr}


\begin{abstract}
Recently, T. Kim considered Euler zeta function which interpolates Euler polynomials at negative integer (see \cite{03}).
In this paper, we study degenerate Euler zeta function which is holomorphic function on complex $s$-plane associated with degenerate Euler polynomials at negative integers.
\end{abstract}

\maketitle

\section{Introduction}

As is well known, the {\it{Euler polynomials}} are defined by the generating function to be

\begin{equation}\label{1}
\frac{2}{e^t+1}e^{xt}=\sum_{n=0} ^{\infty} E_n(x)\frac{t^n}{n!},{\text{ (see [1-7])}}
\end{equation}
When $x=0$, $E_n=E_n(0)$ are called {\it{Euler numbers}}.

For $s\in{\mathbb{C}}$, Kim considered Euler-zeta function which is defined by
\begin{equation}\label{2}
\zeta_E(s,x)=2\sum_{n=0} ^{\infty}\frac{(-1)^n}{(n+x)^s},~(x\neq0,-1,-2,\cdot\cdot\cdot).
\end{equation}
Thus, he obtained the following equation:
\begin{equation}\label{3}
\zeta_E(-n,x)=E_n(x),~(n\in{\mathbb{N}}\cup\left\{0\right\}), \ {\text{ (see \cite{03})}}.
\end{equation}

L. Carlitz introduced {\it{degenerate Euler polynomials}} which are defined by the generating function to be
\begin{equation}\label{4}
\frac{2}{\il^{\iol}+1}\il^{\frac{x}{\lambda}}=\sum_{n=0} ^{\infty} {\mathcal{E}}_n(x|\lambda)\frac{t^n}{n!},{\text{ (see \cite{01})}}.
\end{equation}
When $x=0$, ${\mathcal{E}}_n(\lambda)={\mathcal{E}}_n(0|\lambda)$ are called {\it{degenerate Euler numbers}}.

From \eqref{4}, we note that $\lim_{\lambda\rightarrow0}{\mathcal{E}}_n(x|\lambda)=E_n(x)$, $(n\geq0)$. By \eqref{4}, we easily get
\begin{equation}\label{5}
{\mathcal{E}}_m(\lambda)+{\mathcal{E}}_m(n+1|\lambda)=2\sum_{l=0} ^n(-1)^l(l|\lambda)_m,
\end{equation}
where $(l|\lambda)_m=l(l-\lambda)\cdots(l-\lambda(m-1))$ and $(l|\lambda)_0=1$.

In this paper, we construct degenerate Euler zeta function which interpolates degenerate Euler polynomials at negative integers.
\bigskip
\bigskip

\section{Degenerate Euler zeta function}

For $s\in{\mathbb{C}} \smallsetminus\left\{0,-1,-2,\cdot\cdot\cdot\right\}$, we recall that {\it{gamma function}} is given by
\begin{equation}\label{6}
\Gamma(s)=\ino e^{-t}t^{s-1}dt.
\end{equation}
Let
\begin{equation*}
F(t,x)=\frac{2}{e^t+1}e^{xt}=\sum_{n=0} ^{\infty} E_n(x)\frac{t^n}{n!}.
\end{equation*}
Then, by \eqref{6}, we get
\begin{equation}\label{7}
\begin{split}
\frac{1}{\Gamma(s)}\ino F(-t,x)t^{s-1}dt&=\frac{2}{\Gamma(s)}\ino\frac{1}{1+e^{-t}}e^{-xt}t^{s-1}dt\\
=&\frac{2}{\Gamma(s)}\sum_{m=0} ^{\infty}(-1)^m\ino e^{-(m+x)t}t^{s-1}dt\\
=&\frac{2}{\Gamma(s)}\sum_{m=0} ^{\infty}\frac{(-1)^m}{(m+x)^s}\ino e^{-y}y^{s-1}dy=2\sum_{m=0} ^{\infty}\frac{(-1)^m}{(m+x)^s}\\
=&\zeta_E(s,x).
\end{split}
\end{equation}
From \eqref{7}, we note that
\begin{equation*}
\zeta_E(-n,x)=E_n(x),~(n\in{\mathbb{N}}\cup\left\{0\right\}).
\end{equation*}
Let
\begin{equation}\label{8}
F(t,x|\lambda)=\frac{2}{\il^{\iol}+1}\il^{\frac{x}{\lambda}}=\sum_{n=0} ^{\infty}{\mathcal{E}}_n(x|\lambda)\frac{t^n}{n!}.
\end{equation}
For $\lambda\in(0,1)$ and $s\in{\mathbb{C}}$ with $R(s)>0$, we define degenerate $\Gamma$-function as follows:
\begin{equation}\label{9}
\Gamma(s|\lambda)=\ino \il^{-\iol}t^{s-1}dt.
\end{equation}
Note that $\lim_{\lambda\rightarrow0}\Gamma(s|\lambda)=\Gamma(s)$. From \eqref{2}, we can derive
\begin{equation}\label{10}
\begin{split}
\Gamma(s+1|\lambda)=&\ino\il^{-\iol}t^sdt\\
=&-\frac{s}{\lambda-1}\ino\il^{-\frac{1-\lambda}{\lambda}}t^{s-1}dt\\
=&\frac{s}{1-\lambda}\ino\left(1+\frac{\lambda}{1-\lambda}(1-\lambda)t\right)^{-\frac{1-\lambda}{\lambda}}t^{s-1}dt\\
=&\frac{s}{(1-\lambda)^{s+1}}\ino\left(1+\frac{\lambda}{1-\lambda}y\right)^{-\frac{1-\lambda}{\lambda}}y^{s-1}dy\\
=&\frac{s}{(1-\lambda)^{s+1}}\Gamma\left(s\left|\frac{\lambda}{1-\lambda}\right.\right).
\end{split}
\end{equation}
Therefore, by \eqref{10}, we obtain the following theorem.\\

\bigskip
\begin{theorem}\label{thm1}
For $s\in{\mathbb{C}}$ with $R(s)>0$ and $\lambda\in(0,1)$, we have
\begin{equation*}
\Gamma(s+1|\lambda)=\frac{s}{(1-\lambda)^{s+1}}\Gamma\left(s\left|\frac{\lambda}{1-\lambda}\right.\right).
\end{equation*}
\end{theorem}

\medskip

By Theorem \ref{thm1}, we get
\begin{equation}\label{11}
\begin{split}
\Gamma\left(s\left|\frac{\lambda}{1-\lambda}\right.\right)=&\frac{s-1}{\left(1-\frac{\lambda}{1-\lambda}\right)^s}\Gamma\left(s-1\left|\frac{\frac{\lambda}{1-\lambda}}{1-\frac{\lambda}{1-\lambda}}\right.\right)\\
=&\frac{(s-1)(1-\lambda)^s}{(1-2\lambda)^s}\Gamma\left(s-1\left|\frac{\lambda}{1-2\lambda}\right.\right),
\end{split}
\end{equation}
and
\begin{equation}\label{12}
\Gamma\left(s-1\left|\frac{\lambda}{1-\lambda}\right.\right)=\frac{(s-2)(1-2\lambda)^{s-1}}{(1-3\lambda)^{s-1}}\Gamma\left(s-2\left|\frac{\lambda}{1-3\lambda}\right.\right).
\end{equation}
Continuing this process, we have
\begin{equation}\label{13}
\Gamma\left(s-n\left|\frac{\lambda}{1-\lambda}\right.\right)=\frac{(s-(n+1))(1-(n+1)\lambda)^{s-n}}{(1-(n+2)\lambda)^{s-n}}\Gamma\left(s-(n+1)\left|\frac{\lambda}{1-(n+2)\lambda}\right.\right).
\end{equation}
Thus, by Theorem \ref{thm1}, we get
\begin{equation}\label{14}
\begin{split}
\Gamma(s+1|\lambda)
=&\frac{s}{(1-\lambda)^{s+1}}\Gamma\left(s\left|\frac{\lambda}{1-\lambda}\right.\right)=\frac{s(s-1)}{(1-\lambda)(1-2\lambda)^s}\Gamma\left(s-1\left|\frac{\lambda}{1-2\lambda}\right.\right)\\
=&\frac{s(s-1)(s-2)}{(1-\lambda)(1-2\lambda)(1-3\lambda)^{s-1}}\Gamma\left(s-2\left|\frac{\lambda}{1-3\lambda}\right.\right)=\cdots\\
=&\frac{s(s-1)(s-2)\cdots(s-(n+1))}{(1-\lambda)(1-2\lambda)\cdots(1-(n+1)\lambda)}\left(\frac{1}{1-(n+2)\lambda}\right)^{s-n}\\
&\times\Gamma\left(s-(n+1)\left|\frac{\lambda}{1-(n+2)\lambda}\right.\right).
\end{split}
\end{equation}
Therefore, by \eqref{14}, we obtain the following theorem.

\bigskip

\begin{theorem}\label{thm2}
For $n\in{\mathbb{N}}$, we have
\begin{equation*}
\frac{\Gamma(s+1|\lambda)}{\Gamma\left(s-(n+1)\left|\frac{\lambda}{1-(n+2)\lambda}\right.\right)}=\frac{s(s-1)(s-2)\cdots(s-(n+1))}{(1-\lambda)(1-2\lambda)\cdots(1-(n+1)\lambda)}\left(\frac{1}{1-(n+2)\lambda}\right)^{s-n}.
\end{equation*}
\end{theorem}

\medskip

Let us take $s=n+2$ $(n\in{\mathbb{N}})$ and $\lambda\in\left(0,\frac{1}{n+3}\right)$. Then we have
\begin{equation}\label{15}
\begin{split}
\Gamma(n+3|\lambda)
=&\frac{(n+2)!}{(1-\lambda)(1-2\lambda)\cdots(1-(n+1)\lambda)}\left(\frac{1}{1-(n+2)\lambda}\right)^2\\
&\times\Gamma\left(1\left|\frac{\lambda}{1-(n+2)\lambda}\right.\right).
\end{split}
\end{equation}
For $\lambda\in\left(0,\frac{1}{n+3}\right)$, we observe that
\begin{equation}\label{16}
\begin{split}
&\Gamma\left(1\left|\frac{\lambda}{1-(n+2)\lambda}\right.\right)=\ino\left(1+\left(\frac{\lambda}{1-(n+2)\lambda}\right)t\right)^{-\frac{1-(n+2)\lambda}{\lambda}}dt\\
=&\left.\left(\frac{1-(n+2)\lambda}{\lambda}\right)\left(\frac{\lambda}{(n+3)\lambda-1}\right)\left(1+\left(\frac{\lambda}{1-(n+2)\lambda}\right)t\right)^{\frac{(n+3)\lambda-1}{\lambda}}\right|_0 ^{\infty}\\
=&\left(\frac{1-(n+2)\lambda}{\lambda}\right)\left(\frac{\lambda}{(n+3)\lambda-1}\right)(-1)=\frac{1-(n+2)\lambda}{1-(n+3)\lambda}.
\end{split}
\end{equation}
Thus, by \eqref{15} and \eqref{16}, we get
\begin{equation}\label{17}
\Gamma(n+3|\lambda)=\frac{(n+2)!}{(1-\lambda)(1-2\lambda)\cdots(1-(n+2)\lambda)(1-(n+3)\lambda)}.
\end{equation}
Therefore, by \eqref{17}, we obtain the following theorem.\\

\bigskip

\begin{theorem}\label{thm3}
For $n\in{\mathbb{N}}$, $\lambda\in\left(0,\frac{1}{n}\right)$, we have
\begin{equation*}
\Gamma(n|\lambda)=\frac{(n-1)!}{(1-\lambda)(1-2\lambda)\cdots(1-n\lambda)}.
\end{equation*}
\end{theorem}

\bigskip

In the viewpoint of \eqref{7}, we define {\it{degenerate Euler zeta function}} as follows:
\begin{equation}\label{18}
\zeta_E(s,x|\lambda)=\frac{1}{\Gamma(s|\lambda)}\ino F(-t,x|-\lambda)t^{s-1}dt,
\end{equation}
where $\lambda\in(0,1)$ and $s\in{\mathbb{C}}$ with $R(s)>0$.

From \eqref{18}, we have
\begin{equation}\label{19}
\begin{split}
&\frac{1}{\Gamma(s|\lambda)}\ino F(-t,x|-\lambda)t^{s-1}dt=2\sum_{m=0} ^{\infty}(-1)^m\ino\il^{-\frac{m+x}{\lambda}}t^{s-1}dt\\
\ \ \ \ =&\frac{2}{\Gamma(s|\lambda)}\sum_{m=0} ^{\infty}(-1)^m\ino\left(1+\frac{\lambda}{m+x}(m+x)t\right)^{-\frac{m+x}\lambda}t^{s-1}dt\\
\ \ \ \ =&\frac{2}{\Gamma(s|\lambda)}\sum_{m=0} ^{\infty}\frac{(-1)^m}{(m+x)^s}\ino\left(1+\frac{\lambda}{m+x}y\right)^{-\frac{m+x}{\lambda}}y^{s-1}dy\\
\ \ \ \ =&2\sum_{m=0} ^{\infty}\frac{(-1)^m}{(m+x)^s}\frac{\Gamma\left(s\left|\frac{\lambda}{m+x}\right.\right)}{\Gamma(s|\lambda)},
\end{split}
\end{equation}
where $x\neq0,-1,-2,\cdot\cdot\cdot$.

Therefore, by \eqref{18} and \eqref{19}, we obtain the following theorem.

\bigskip

\begin{theorem}\label{thm4}
For $s\in{\mathbb{C}}$ with $R(s)>0$, we have
\begin{equation*}
\zeta_E(s,x|\lambda)=2\sum_{m=0} ^{\infty}\frac{(-1)^m}{(m+x)^s}\frac{\Gamma\left(s\left|\frac{\lambda}{m+x}\right.\right)}{\Gamma(s|\lambda)},
\end{equation*}
where $x\neq0,-1,-2,\cdot\cdot\cdot.$
\end{theorem}

\bigskip

Let $s=n\in{\mathbb{N}}$ and $\lambda\in\left(0,\frac{1}{n}\right)$. Then, we have
\begin{equation}\label{20}
\zeta_E(n,x|\lambda)=2\sum_{m=0} ^{\infty}\frac{(-1)^m}{(m+x)^n}\frac{\Gamma\left(n\left|\frac{\lambda}{m+x}\right.\right)}{\Gamma(n|\lambda)}.
\end{equation}
From Theorem \ref{thm3}, we note that
\begin{equation}\label{21}
\begin{split}
&\frac{\Gamma\left(n\left|\frac{\lambda}{m+x}\right.\right)}{\Gamma(n|\lambda)}\\
=&\frac{(1-\lambda)(1-2\lambda)\cdots(1-n\lambda)}{(n-1)!}\frac{(n-1)!}{\left(1-\frac{\lambda}{m+x}\right)\left(1-\frac{2\lambda}{m+x}\right)\cdots\left(1-\frac{n\lambda}{m+x}\right)}\\
=&\frac{(1-\lambda)(1-2\lambda)\cdots(1-n\lambda)}{\left(1-\frac{\lambda}{m+x}\right)\left(1-\frac{2\lambda}{m+x}\right)\cdots\left(1-\frac{n\lambda}{m+x}\right)}.
\end{split}
\end{equation}
Therefore, by \eqref{20} and \eqref{21}, we obtain the following theorem.\\

\bigskip

\begin{theorem}\label{thm5}
For $n\in{\mathbb{N}}$ and $\lambda\in\left(0,\frac{1}{n}\right)$, we have
\begin{equation*}
\zeta_E(n,x|\lambda)=2\sum_{m=0} ^{\infty}\frac{(-1)^m}{(m+x)^n}\frac{(1-\lambda)(1-2\lambda)\cdots(1-n\lambda)}{\left(1-\frac{\lambda}{m+x}\right)\left(1-\frac{2\lambda}{m+x}\right)\cdots\left(1-\frac{n\lambda}{m+x}\right)}.
\end{equation*}
\end{theorem}

\bigskip

From \eqref{8}, we have
\begin{equation}\label{22}
\begin{split}
&\frac{1}{\Gamma(s|\lambda)}\ino F(-t,x|-\lambda)t^{s-1}dt\\
=&\frac{1}{\Gamma(s|\lambda)}\ino\frac{2}{\il^{-\iol}+1}\il^{-\frac{x}{\lambda}}t^{s-1}dt\\
=&\frac{1}{\Gamma(s|\lambda)}\sum_{m=0} ^{\infty}{\mathcal{E}}_m(x|-\lambda)\frac{(-1)^m}{m!}\ino t^{s+m-1}dt.
\end{split}
\end{equation}
Thus, by \eqref{22}, we get
\begin{equation}\label{23}
\zeta_E(-n,x|\lambda)=\frac{2\pi i}{n!}{\mathcal{E}}_n(x|-\lambda)(-1)^n\frac{1}{\Gamma(-n|\lambda)}.
\end{equation}
Now, we observe that
\begin{equation}\label{24}
\begin{split}
\Gamma(-n|\lambda)=&\ino\il^{-\iol}t^{-n-1}dt\\
=&\frac{2\pi i}{n!}\frac{1}{\lambda}\left(1+\frac{1}{\lambda}\right)\left(2+\frac{1}{\lambda}\right)\cdots\left((n-1)+\frac{1}{\lambda}\right)(-1)^n\lambda^n\\
=&\frac{2\pi i}{n!}(\lambda +1)(2\lambda+1)\cdots((n-1)\lambda+1)(-1)^n.
\end{split}
\end{equation}
From \eqref{23} and \eqref{24}, we have
\begin{equation}\label{25}
\zeta_E(-n,x|\lambda)=\frac{{\mathcal{E}}_n(x|-\lambda)}{(\lambda+1)(2\lambda+1)\cdots((n-1)\lambda+1)}.
\end{equation}
Therefore, by \eqref{25}, we obtain the following theorem.

\bigskip

\begin{theorem}\label{thm6}
For $n\in{\mathbb{N}}\cup\left\{0\right\}$, we have
\begin{equation*}
\zeta_E(-n,x|\lambda)={\mathcal{E}}_n(x|-\lambda).
\end{equation*}
\end{theorem}

\bigskip

\noindent{\bf{Remark.}} We note that $\zeta_E(s,x|\lambda)$ is analytic function in whole complex $s$-plane.

\vskip .1in

\noindent{\bf{Remark.}}
\begin{equation*}
\lim_{\lambda\rightarrow0}\zeta_E(-n,x|\lambda)=E_n(x)=\zeta_E(-n,x).
\end{equation*}

\end{document}